\documentclass{amsart}
\usepackage{amsmath,amssymb,enumerate}
\usepackage{epsf,epsfig,amsfonts,graphicx}  
\usepackage[usenames,dvipsnames]{color}
 \usepackage{hyperref}
\usepackage[T1]{fontenc}  
\usepackage[utf8]{inputenc}

 \date{May 12, 2016}   
  
 \numberwithin{equation}{section}   
\newcommand{\ds}{\displaystyle}

\renewcommand{\r}{\mathbb{R}}

\newtheorem*{maintheoremA}{\rm\bf Theorem A} 
\newtheorem*{maintheoremB}{\rm\bf Theorem B}
\newtheorem*{corollaryC}{\rm\bf Corollary C}
\newtheorem*{theoremD}{\rm\bf Theorem D}
\newtheorem*{theoremBB}{\rm\bf Theorem B'}
\newtheorem*{corollaryCC}{\rm\bf Corollary C'}

\newtheorem{theorem}{\rm\bf Theorem}[section]
\newtheorem{proposition}[theorem]{\rm\bf Proposition}
\newtheorem{lemma}[theorem]{\rm\bf Lemma}
\newtheorem{corollary}[theorem]{\rm\bf Corollary}

\theoremstyle{definition}

\newcommand{\weg}[1]{}

\title{The  Myers-Steenrod theorem for Finsler manifolds of low regularity.}

\author{Vladimir S. Matveev} 
\address{Institut f\"ur Mathematik, Friedrich-Schiller Universit\"at Jena\\
07737 Jena, Germany}  
\email{vladimir.matveev@uni-jena.de}

\author{Marc Troyanov} 
\address{Section de Math{\'e}matiques,  
\'Ecole Polytechnique F{\'e}derale de Lausanne, station 8,
1015 Lausanne - Switzerland} 
\email{marc.troyanov@epfl.ch}

\begin{document}

\bigskip

\begin{abstract}{We prove a version of Myers-Steenrod's theorem for Finsler manifolds
under minimal regularity hypothesis. In particular we  show that an isometry between 
$C^{k,\alpha}$-smooth (or partially smooth) Finsler metrics,  with $k+\alpha>0$,  $k\in \mathbb{N} \cup \{0\}$,
and $0 \leq \alpha \leq 1$   is  necessary a diffeomorphism of class $C^{k+1,\alpha}$.
A generalisation of this result to the case of Finsler 1-quasiconformal mapping is given. 
The proofs are  based on the reduction of the Finlserian problems  to  Riemannian ones  
with the help of the  Binet-Legendre metric.  
 }
  \\ 
\medskip

\noindent 2000 AMS Mathematics Subject Classification:  53b40,53c60,35b65. \\
Keywords:  Finsler metric, isometries,  Myers-Steenrod theorem, Binet-Legendre metric.
\end{abstract}
\thanks{We thank the Friedrich-Schiller-Universit\"at Jena, EPFL and the Swiss National Science Foundation for their support.}
\maketitle


\section{Introduction}

The main goal in this paper is to prove a Myers-Steenrod theorem for Finsler
manifolds under low regularity assumptions. More precisely, we give an answer to
the following question: 
\emph{What is the regularity (smoothness) of a distance preserving  bijection  $\phi : M_1\to M_2$ 
between two  Finsler manifolds 
 $(M_1,F_1)$ and $(M_2,F_2)$ of low regularity?}

\smallskip

This question has been investigated by several authors starting with the seminal 1939 paper by S. Myers and N.  Steenrod
on Riemannian isometries \cite{MyersSteenrod1939}.  We provide a brief historical account at the end of the paper.  
Our approach in this paper is to reduce the Finsler case to the Riemannian one and to use results on Riemannian isometries
proved in   \cite{Lytchak2006,reshetnyak,sabitov1993,shefel,Taylor2006}. Our tool is the \emph{Binet-Legendre} metric $g_{_F}$. 
It is a  Riemannian metric    on $M$ canonically  associated to  $F$ 
by some averaging procedure, see    Section \ref{BL} below. This metric enjoys some natural geometric properties
that makes it a  useful tool in Finsler geometry; this is concretely illustrated in \cite{MT2012,MT2015}. 
We state our main result as the following two Theorems; in these statements, we assume 
\begin{equation}
 k\in \mathbb{N}\cup \{0\},  \quad  0 \leq\alpha\leq 1 \quad \text{and} \quad k+\alpha > 0.
\end{equation}

\smallskip

\begin{maintheoremA}
 Let  $\phi : M_1\to M_2$ be a distance preserving  bijective map   between two  Finsler manifolds  $(M_1,F_1)$ and $(M_2,F_2)$.
 Assume that the Binet-Legendre metrics $g_{_{F_1}}$ and $g_{_{F_2}}$  associated to $F_1$ and $F_2$  are locally  of class 
 $C^{k,\alpha}$. Then $\phi \in C_{loc}^{k+1,\alpha}(M_1,M_2)$ and $\phi^*(F_2) = F_1$.
\end{maintheoremA}

\medskip

Recall that the notation $C^{k, \alpha}$ stands for the class of functions or mappings that are   $k$ times continuous differentiable
and all of  whose partial derivatives of order $k$ are   H\"older continuous of order $\alpha$ (or Lipschitz continuous if $\alpha = 1$).

\begin{maintheoremB} 
   If a Finsler metric $F$ is  of class  $C_{loc}^{k,\alpha}$, then its  Binet-Legendre metric is also of class  $C_{loc}^{k,\alpha}$.
 \end{maintheoremB}

Combining these two theorems, we obtain the following: 
\begin{corollaryC} \label{corC} 
Any distance preserving bijective map  $\phi$  between two manifolds with H\"older continuous Finsler metrics
is  a diffeomorphism. Moreover, if the metrics are of class  $C_{loc}^{k,\alpha}$ with $k+\alpha>0$, then the  diffeomorphism
$\phi$   is  of class $C_{loc}^{k+1,\alpha}$. 
\end{corollaryC}

\medskip

It is important to note that  the Binet-Legendre metric $g_{_F}$ of a Finsler metric $F$  may be smooth even if the metric $F$ is of low
regularity. A simple example is a  Minkowski metric on $\r^n$ with polyhedral unit ball. This Finsler metric is not of class $C^1$, but its associated Binet-Legendre metric is  an Euclidean metric and is therefore $C^{\infty}$. 
Other examples are given by the Funk metrics and, more generally, the so called Zermelo metrics,  see  \cite[\S 5.2]{MT2012}. 
Additional examples are given by the class of partially smooth Finsler metrics introduced in  \cite[\S 2]{MT2015},
see also  \S \ref{rec.reg} below.

\medskip

The paper is organized as follows. In \S \ref{sec2}, we first recall some basic facts on Finsler  structures  and we give 
the definition of  Binet-Legendre metric and its basic properties. We then state the optimal regularity result for Riemannian 
isometries    due to   Yu.  Reshetnyak, I. Sabitov,  M. Taylor, A. Lytchak and A. Yaman. In \S \ref{sec3} we prove Theorem A
and we discuss a generalization to $1$-quasiconformal maps between Finsler manifolds. In \S  \ref{sec4} we recall   some  fine properties
of H\"older continuous maps and in \S \ref{rec.reg} we use them to  prove Theorem B and  give some sufficient conditions for the regularity of the
Binet-Legendre metric.  In \S \ref{secAppl} we give an application to Finsler manifolds admitting non trivial dilations  and  discuss the group of isometries of a Finsler manifold.
The final section contains a brief history of the Myers-Steenrod Theorem  for Riemannian manifolds  and references on the subject.

 \medskip

We would like to thank CY Guo, S. Ivanov, A. Lytchak , D.  Repovš, I. Kh. Sabitov, and E.  Ščepin, for useful discussions. 
We are particularly grateful to A. Lytchak for attracting our interest to the paper \cite{sabitov1993}.

\section{Preliminaries}
 \label{sec2} 
\subsection{Finsler metrics} 

A  \emph{Finsler structure} on a domain\footnote{Finsler structure are more generally defined as  continuous maps
$F$ on the tangent bundle of a $C^1$ manifold $M$ satisfying the conditions (a), (b) and (c)  
(for a domain $ \mathcal{U} \subset \r^n$, we identify $T\mathcal{U}$  with $\mathcal{U}\times \r^n$).  To define a $C^{k,\alpha}$ Finlser metric
on a  manifold, we need to  assume that an atlas of class $C^2 \cap C^{k+1,\alpha}$ is given on the manifold.} $\mathcal{U} \subset \r^n$ is a continuous 
function  $F : \mathcal{U}\times \r^n \to [0,\infty)$
such that for any $x\in \mathcal{U}$ and any $v, w \in \r^n$ we have
\begin{enumerate}[ \ (a)]
  \item $F(x,\lambda \cdot v) = \lambda \cdot   F (x,v) $ for any $\lambda \geq 0$.
   \item $F (x,v+ w ) \le F(x,v) + F(x,w)$.
   \item $F(x,v)= 0 $ \  $ \Rightarrow$ \  $v=0$.  
\end{enumerate} 

\medskip
Before proceeding, let us  recall a few more definitions:

\smallskip \newpage 

\textbf{Definitions} 
\begin{enumerate}[i)]
  \item The Finsler structure $F$ on $\mathcal{U}$ is said to be \emph{reversible} if $F(x, -v) = F(x,v)$
  for any $(x,v) \in  \mathcal{U}\times \r^n$.
  \item If the Finsler structure $F$ is independent from the point $x\in \mathcal{U}$, one says that it is a \emph{Minkowski norm}.
  \item The Finsler structure is said to be \emph{of class} $C^{k, \alpha}$ if  the restriction of $F$ to an open neighborhood of
  $\mathcal{U}\times S^{n-1} = \{(x,v) \in \mathcal{U} \times \r^n  \mid |v| = 1\}$  is a function of class  $C^{k, \alpha}$.
\end{enumerate}

\medskip

\textbf{Remarks}  (i) Observe   that a Minkowski norm on $\r^n$  is a norm in the usual sense if and only if it is reversible. \\
(ii) It is often assumed in the Finsler literature (but not in the present paper) that 
 the Finsler structure  is of class  at least $C^2$  and for any $(x,v)\in \mathcal{U}\times (\r^n\setminus \{0\})$   the vertical Hessian matrix
\begin{equation*}\label{vert.hessian}
   \left(\frac{\partial^2 F^2}{\partial v_i \partial v_j}  \right)
\end{equation*}
 is  positive definite. Such metrics are called \emph{strongly convex}.
 Note that most results that have been proved so far on the  smoothness of Finlser isometries
 assume the metrics to be  strongly convex, see \S \ref{secHistory}.

\medskip

Given a Finsler structure $F$ on $\mathcal{U} \in \mathbb{R}^n$, one defines the $F$-length of a $C^1$-curve
$\gamma : [a,b] \to \mathcal{U}$ as
$$
  \ell_{F} (\gamma) = \int_a^b F(\gamma (t), \dot\gamma(t)) dt.
$$
The \emph{length} of any curve is invariant under sense preserving reparametrisation, it is 
also invariant under sense reversing reparametrisation if the Finsler structure is reversible.
The \emph{distance} between two points $x$ and $y$ in $\mathcal{U}$ is defined as the infimum of the
length of all $C^1$ curves joining them
$$
  d_F(x,y) = \inf \{\ell_{F} (\gamma) \mid \gamma \in C^1([0,1], \mathcal{U}), \gamma(0) = x, \gamma(1) = y\}.
$$
The distance satisfies $d_F(x,y) + d_F(y,z) \geq d_F(x,z)$ and $d_F(x,y) =0$ if and only if $x=y$.
If the Finsler structure is reversible, then the distance is also symmetric:  $d_F(x,y) =d_F(y,x)$.

\begin{lemma}\label{lemma.bilipshitz}
Let $\mathcal{U}$ be a convex domain in $\r^n$ and $F$ be a $C^0$ Finsler metric on $\mathcal{U}$. 
Suppose that there exists a constant $C>0$ such that 
$$
  \frac{1}{C}|v | \leq F(x,v) \leq  C |v|
$$
for all $(x, v) \in \mathcal{U}\times \r^n$,  where $|\cdot|$ denotes the usual Euclidean norm. 
Then the Finsler distance is bilipschitz equivalent to the 
Euclidean distance in $\mathcal{U}$, more precisely we have
\begin{equation}\label{eq.bbilip}
   \frac{1}{C}|q-p | \leq d_F(p,q)\leq  C |q-p|
\end{equation}
for any $p,q \in \mathcal{U}$.
\end{lemma}

\textbf{Proof.}   Fix two points $p$ and $q$ in $\mathcal{U}$ and  choose a $C^1$
path $\gamma : [0,1] \to \mathcal{U}$ joining them. We then have
$$
  |q-p| \leq \int_0^1 |\dot \gamma(t)| dt \leq C \int_0^1 F(\gamma(t), \dot \gamma(t)) dt.
$$
Taking the infimum of this inequality on all paths from $p$ to $q$ one obtains 
$|q-p| \leq d_F(p,q)$.  To prove the other inequality, we consider the segment
$\alpha(t) = t q + (1-t)p$, this path is contained in $\mathcal{U}$ since this is a convex domain.
We then have
$$ d_F(p,q) \leq  \int_0^1 F(\alpha(t), \dot \alpha(t)) dt \leq C  \int_0^1 |\dot \alpha(t)| dt 
 = C  \int_0^1 |q-p| dt  = C |q-p|.
$$

\qed

\textbf{Remark.}   This lemma also holds for a quasi-convex domain but with a different constant
in (\ref{eq.bbilip}). Recall that a domain $\mathcal{U}$ in $\r^n$ is   quasi-convex if there is a constant $K$
such that for any $x,y \in \mathcal{U}$ there exists a path in $\mathcal{U}$ joining $x$ to $y$ of Euclidean length at
most $K |y-x|$.

\subsection{The Binet-Legendre metric associated to a Finsler structure} \label{BL}

Given a Finsler structure $F$ on $\mathcal{U}$, one defines for any point $x \in \mathcal{U}$ 
the \emph{associated Finsler unit tangent ball}  
$$
  \Omega_x=  \{ v \in \r^n \mid  F(x,v) < 1\}.
$$
The boundary of $\Omega_x$ is often called the \emph{indicatrix} of $F$ at $x$.

\medskip

\textbf{Definition}  The \emph{Binet-Legendre} metric associated to the Finsler structure $F$ on $\mathcal{U} \subset \r^n$
is the  Riemannian metric whose metric tensor $g_{ij}$ is  defined at any point $x \in \mathcal{U}$ to be the inverse
matrix of
\begin{equation}\label{eq.BL1}
  g^{ij}(x) = \frac{(n+2)}{\lambda^n(\Omega_x)} \int_{\Omega_x} v_iv_j \; dv_1\dots dv_n,
\end{equation}
where $\lambda^n(\Omega_x)$ is the  Lebesgue measure of $\Omega_x \subset \r^n$. 

\medskip

We refer to \cite{MT2012} for an intrinsic (coordinate free) definition of  the Binet-Legendre metric, discussions of its properties and various applications to Finsler geometry. This metric first appeared in the work of Paul Centore \cite{Centore}, who called it the \emph{osculating metric}.

\medskip

It will be useful to rewrite formula (\ref{eq.BL1}) in polar  coordinates. Let us write 
$v = r u$ with $u \in S^{n-1}$ and $r \in \r_+$.  We then have
$$
 \lambda^n(\Omega_x) =   \int_{\Omega_x} dv = 
  \int_{S^{n-1}} \left( \int_0^{1/{F(x,u)}} r^{n-1} dr   \right) d\sigma(u)
 = \frac{1}{n} \int_{S^{n-1}} \frac{d\sigma(u)}{F(x,u)^n}, 
$$
where $d\sigma$ is the standard volume form on $S^{n-1}$.
this proves that $x \mapsto \lambda^n(\Omega_x)$ is of class  $C^{k, \alpha}$. 

\smallskip

We likewise have 
\begin{align*}
\int_{\Omega_x} v_iv_j \; dv &=  \int_{S^{n-1}} \left( \int_0^{1/{F(x,u)}}(r^2u_iu_j) r^{n-1} dr   \right) d\sigma(u)
     \\&= 
    \int_{S^{n-1}} u_iu_j\left( \int_0^{1/{F(x,u)}} r^{n+1} dr   \right) d\sigma(u)
    \\&=  \frac{1}{(n+2)} \int_{S^{n-1}} \frac{u_iu_j}{F(x,u)^{n+2}} \,  d\sigma(u).
\end{align*}

\smallskip

It then follows that (\ref{eq.BL1}) can be written as
\begin{equation}\label{eq.BL2}
    g^{ij}(x) =   \frac{n {\ds \int_{S^{n-1}}}  {u_iu_j}{F(x,u)^{-(n+2)}}d\sigma}{ {\ds \int_{S^{n-1}}} F(x,u)^{-n}d\sigma}.
\end{equation}

\medskip

Some basic  properties of the Binet-Legendre metric are stated and proved in  \cite{MT2012}.
Let us recall  that if the Finsler metric $F$ on $\mathcal{U}$ is Riemannian, that is if there exists a Riemannian
metric $h$ on $\mathcal{U}$ such that $F(x,v) = \sqrt{h_x(v,v)}$ for all $(x,v)\in \mathcal{U}\times \r^n$, 
then the Binet-Legendre metric of $F$ coincide with $h$.

\subsection{Regularity of isometries in the Riemannian case}
We will use the following result on Riemannian manifolds:

 \begin{theorem}\label{th.regRiem}
 Let $(\mathcal{U}_1,g_1)$ and $(\mathcal{U}_2, g_2)$ be  two domains of $\r^n$  equipped with   
 Riemannian metrics  $g_1$ and $g_2$  of class $C^{k, \alpha}$ with $k\in \mathbb{N}\cup \{0\}$, $0\leq \alpha  \le  1$
 and $k+\alpha>0$.
  Let $\phi : (\mathcal{U}_1,g_1) \to (\mathcal{U}_2, g_2)$ be a  bijective mapping.  Then, the following  
  conditions are equivalent:
\begin{enumerate} \item $\phi$ is a distance-preserving homeomorphism,
\item  $\phi$  is bi-Lipschitz and  $\phi^* {g_2}= {g_1}$ almost everywhere, \label{item2} 
\item  $\phi$  is a diffeomorphism of class $C^{r+1, \alpha}$ and $\phi^* {g_2}= {g_1}$. \label{item3}
  \end{enumerate} 
\end{theorem}

\medskip

This theorem is the combined result of the work of several authors,   
we explain this now and  give a historical overview in section \ref{secHistory} below.

\medskip
 
\textbf{Sketch of proof.} 
The implication (1) $\Rightarrow$ (2)  holds in fact even for $C^0$ metrics $g_1$ and $g_2$
as explained in the first step of the proof of  \cite[Theorem 2.1]{Taylor2006},  and (3) $\Rightarrow$ (1) 
is obvious. The main issue  is therefore the implication  (2) $\Rightarrow$ (3).

\smallskip

If  $k=0$ and $0 < \alpha < 1$, then the implication  (2) $\Rightarrow$ (3) has been proved by Yu. Reshetnyack in \cite[Theorem 2]{reshetnyak} (where in fact he proves a more general statement). 
An independant proof is given by M. Taylor in  \cite{Taylor2006}.

\smallskip

If  $k+\alpha = 1$, that is either  $k=0$ and $\alpha = 1$, or $k=1$ and $\alpha = 0$, then the same implication is
proved   by  I. Sabitov in  \cite{sabitov1993}. 

\smallskip

The general case $k\geq 1$ and $0 \leq \alpha \leq 1$ follows now from the following argument of Calabi and Hartman \cite[\S 5]{Calabi1970}. Since $k\geq 1$, the map $\phi$ is  $C^2$ and satisfy $\phi^*g_2 = g_1$ from 
the previous case. 
Recall now that in any local coordinates systems, the Christoffel symbols
 $\tilde{\Gamma}^{m}_{\nu \lambda}$ and  $\tilde{\Gamma}^{m}_{\nu \lambda}$ of the metrics $g_1$ and $g_2$ are related by the formula
\begin{equation} \label{eq.ChangeCS}
 \frac{\partial^2\phi^m}{\partial x^i \partial x^j} = \Gamma^{\mu}_{ij} \frac{\partial \phi^m}{\partial x^{\mu}}
 - \tilde{\Gamma}^{m}_{\nu \lambda} \frac{\partial \phi^{\nu}}{\partial x^{i}} \frac{\partial \phi^{\lambda}}{\partial x^{j}}.
\end{equation}
Since the Christoffel symbols are given by algebraic expressions involving the first derivatives of the metric
tensors $g_{1}$ and $g_{2}$, the implication  (2) $\Rightarrow$ (3) follows from (\ref{eq.ChangeCS})
by induction on $k$ for any given $0 \leq \alpha \leq 1$.

\qed

\medskip
   
\textbf{Remarks \  1.)} This Theorem is an optimal regularity result for isometries of  $C^{k, \alpha}$  Riemannian metrics,  which of course implies that our Theorem A, and hence Corollary C, are also optimal. Indeed, E.  Calabi and P. Hartman have built an example of a non differentiable isometry  between two $C^0$ Riemannian metrics, see \cite[\S 6]{Calabi1970}.

\smallskip

\textbf{2.)} The result is also optimal in the sense that one cannot expect an  isometry between $C^{k,\alpha}$ metrics to be more than  $C^{k+1,\alpha}$-regular. A simple (Riemannian) example can be build as follows:  choose a diffeomorphism $f : \r \to \r$ which is of class $C^{k+1,\alpha}$
but does not belong to  $C^{k+1,\beta}$ for any $\beta > \alpha$. Consider the following Riemannian metrics on the plane: 
$$g_1= f'(x_1)^2 dx_1^2 + dx_2^2 \quad  \textrm{and} \quad  g_2=dx_1^2 +dx_2^2.$$  
Then the map $\phi:(x_1,x_2)\mapsto  (f(x_1),x_2)$ is an isometry from $(\r^2, g_1)$ to $(\r^2,g_2)$. By construction, the map $\phi$  is of class  $C^{k+1,\alpha}$  but not of class $C^{k+1,\beta}$
for any $\beta>\alpha$.

\smallskip

\textbf{3.)} Note that the implication  (2) $\Rightarrow$ (3) fails  if we only assume the map $\phi$ to be differentiable almost everywhere.
Consider for instance\footnote{We owe this example to Changyu Guo.}  the map from the unit square $(0,1)^2$ to $\r^2$ defined by $\phi(x,y) = (x,y+\nu (x))$, where $\nu : (0,1) \to (0,1)$ is the Cantor-Lebesgue function. Then $\phi$ is a homeomorphism onto its image such that $\phi^* g_0 = g_0$ a.e. where $g_0 = dx^2+dy^2$ is the standard Euclidean metric, yet $\phi \not\in C^1$. However, Reshetnyak's theorem tells us that assuming $\phi \in W^{1,n}$ is sufficient for the implication (2) $\Rightarrow$ (3).

\section{Proof of Theorem A} \label{mainTh} \label{sec3}
 {
 \medskip

We first observe that our hypothesis, together with Lemma \ref{lemma.bilipshitz},  imply that 
the map $\phi$ is bilipschitz. It then follows from Rademacher's theorem that $\phi$ is 
(Frechet) differentiable almost everywhere  \cite[\S 3.1.2]{EvansGariepy}. 

\smallskip

We  claim that \emph{the following condition holds at any point $x\in \mathcal{U}_1$ where $\phi$
is differentiable:}
$$
  F_1(x, v) = F_2(\phi(x), d\phi_x(v)), \qquad  (\forall v \in \r^n).
$$
To prove the claim, we may suppose that $\mathcal{U}_1$ and $ \mathcal{U}_2$ contain the origin and that $\phi(0) = 0$. 
Suppose also that $\phi$ is differentiable at $x=0$ and let us introduce the following notations (for $t>0$ and $i = 1,2$)
$$
 \mathcal{U}_{i,t} = \frac{1}{t} \mathcal{U}_i,  \quad   d_{i,t}(x,y) = \frac{1}{t} d_i(tx, ty),   \quad  F_{i,t}(x,v) =  F_{i}(tx,v) =  \frac{1}{t} F_{i}(tx,tv).
$$
Let us also define the map $\phi_t :  \mathcal{U}_{1,t} \to  \mathcal{U}_{2,t}$ by $\phi_t(x)   = \frac{1}{t} \phi(tx)$.
From the definitions, we see that  $d_{i,t}(x,y)$ is the distance on $\mathcal{U}_{i,t}$  associated to the Finsler metric $F_{i,t}(x,v)$. 
It is also easy to check that 
$\phi_t$ is a distance preserving bijection from $(\mathcal{U}_{1,t},d_{1,t})$ to $(\mathcal{U}_{2,t},d_{2,t})$.

\smallskip

Since we assumed  that $\phi$ is  differentiable at $0$, the limit $\phi_{0} = \lim_{t\to 0} \phi_t$ exists uniformly, and 
it is a map from $\r^n$ to $\r^n$. In fact we have $\phi_{0} = d\phi_0$ (the differential of $\phi$ at $0$), in particular it is a linear map.

\smallskip

By continuity  $\phi_{0}$ is an isometry for the limit distances  $d_{1,0} = \lim_{t\to 0}d_{1,t}$ 
and  $d_{2,0} = \lim_{t\to 0}d_{2,t}$.
These distances are associated to the Finsler metrics  $F_{i,0}(x,v) = \lim_{t\to 0} F_{i,t}(x,v) = F_{i}(0,v)$. 
Note that in fact $F_{1,0}$ and  $F_{2,0}$ are  Minkowski metric on  $\r^n$.
This implies  that $d\phi_0=\phi_0$ is an isometry between two Minkowski spaces $(\r^n, F_{1,0})$ and  $(\r^n, F_{2,0})$.
The claim is proved.

\smallskip

Let us denote by $g_i$ the Binet-Legendre metric of $F_i$. By hypothesis, this  is a $C^{k, \alpha}$ Riemannian metric  on 
$\mathcal{U}_i$. It follows now from  the claim and standard properties of the Binet-Legendre metric that  
$$
  \left.g_2\right|_{\phi(x)} (d\phi (v),d\phi(w)) =   \left.g_1\right|_x(v,w),  
$$
for  a.e. $x\in \mathcal{U}_1$ and all $v,w \in \r^n$.  
\smallskip

We thus have proved that $g_i$ is a $C^{k, \alpha}$ Riemannian metric on $\mathcal{U}_i$ for $i=1,2$ and 
$\phi : \mathcal{U}_1 \to \mathcal{U}_2$ is a bilipschitz map such that $\phi^*g_2 = g_1$ almost everywhere. 
By Theorem \ref{th.regRiem}  we conclude that $\phi \in C^{k+1, \alpha}(\mathcal{U}_1,\mathcal{U}_2)$. 
Theorem A is proved.

\qed

\subsection{Extension of Theorem A to the case of 1-quasiconformal maps}

In this section we briefly explain how to extend our results to the case of
quasiconformal maps. Recall that a homeomorphism $\phi : X \rightarrow Y$
between two metric measure spaces $(X, d)$ and $(Y, d')$ is said to be
\emph{quasiconformal} if the \emph{linear distorsion}, defined as
$$
H ( x, \phi) = \limsup_{r \rightarrow 0} \frac{\max\{d'(\phi (x),\phi (z)) \mid z \in B (x,r) \}}{\min\{ d'(\phi (x),\phi (z)) \mid
   z \in B (x,r)\}} 
$$  
is uniformly bounded in $X$. The map is said to be  \emph{1-quasiconformal} if $H (x, \phi) = 1$ a.e. 
A regularity theorem  for Riemannian $1$-quasiconformal has been proved by  Reshetnyak in  \cite{reshetnyak}. See also 
\cite{Iwaniec1982} and  \cite{LS}. Using this result we can state the following 

\begin{theoremD}
 Let  $\phi : (M_1,F_1)\to (M_2,F_2)$ be a $1$-quasiconformal map   between two  Finsler manifolds.
 Assume that the Binet-Legendre metrics $g_1$ and $g_2$  of $F_1$ and $F_2$  belong to 
 $C_{loc}^{k,\alpha}$, where $k\in \mathbb{N}\cup \{0\}$, \  $0 \leq\alpha < 1$ and $k+\alpha > 0$.
Then $\phi \in C_{loc}^{k+1,\alpha}(M_1,M_2)$ and $\phi^*(F_2) = F_1$.
In particular, the map $\phi$ is conformal.
\end{theoremD}

\medskip

\textbf{Sketch of proof.} \ The notion of quasiconformality is invariant
under a bilipschitz change of the metrics, therefore the restriction of the
quasi-conformal map $\phi : M_1 \rightarrow M_2 $ to any coordinate
neighborhood is also quasiconformal with respect to the Euclidean metric.
Using the standard theory of quasiconformal maps in $\mathbb{R}^n$, one then
deduces that $\phi$ is differentiable almost everywhere and that its
differential is a.e. invertible (see e.g.  \cite[Chapter 6]{IwaniecMartin2001}, the
differentiability a.e. follows from Rademacher-Stepanov's theorem and the fact
that the differential is a.e. invertible follows form the change of variable
formula in integrals).
Now using a blow-up argument as in the beginning of the proof of Theorem A, \
we obtain
\[ \max_{F_1 ( x, u) = 1} F_2 ( \phi ( x), d \phi_x ( u)) = H ( x, \phi)
   \cdot \min_{F_1 ( x, u) = 1} F_2 ( \phi ( x), d \phi_x ( u)) \]
at any point $x \in M_1$ where $\phi$ is differentiable with invertible
differential. In particular, if $\phi : M_1 \rightarrow M_2 $ is
1-quasiconformal, then the following equality holds for almost all $x \in M_1$
and all $v \in T_x M$
\[ F_2 ( \phi ( x), d \phi_x ( v)) = \mu ( x) F ( x, v), \]
where 
$$\mu ( x) = \max_{F_1 ( x, u) = 1} F_2 ( \phi ( x), d \phi_x ( u)) =
\min_{F_1 ( x, u) = 1} F_2 ( \phi ( x), d \phi_x ( u)).$$
The above identity can also be written as $\phi^{\ast} F_2 = \mu F_1$ a.e.  This implies that
$\phi^{\ast} g_2 = \mu g_1$ a.e., that is \ $\phi : ( M_1, g_1) \rightarrow (M_2, g_2) $ 
is Riemannian $1$-quasiconformal. We conclude from \cite[Theorem 2]{reshetnyak}
together with \cite[Theorem 4.5]{LS} or  \cite{Iwaniec1982,NikolaevShefel1986} that $\phi \in C_{loc}^{k + 1, \alpha}$.

\qed

\section{On H\"older  and Lipschitz maps} \label{sec4} 

In this section we collect some  facts on H\"older continuous and Lipschitz mappings
that we will need later.
Let $\mathcal{U}$ be a domain in $\r^n$ and $0 < \alpha <1$. A map $f : \mathcal{U} \to \r^m$ 
is said to be \emph{H\"older continuous} of class $C^{0,\alpha}$ if 
\begin{equation}\label{HolderSemiNorm}
  [f]_{C^{0,\alpha}(\mathcal{U})} =
   \sup \left\{ \  \frac{|f(y)-f(x)|}{|y-x|^{\alpha}}  \,  \big| \,  x,y \in \mathcal{U}, \ x\neq y  \right\}  < \infty.
\end{equation}
Here $|y-x|$ denotes the standard Euclidean distance between $x$ and $y$.
The map $f$ is said to be \emph{Lipschitz} if the same condition holds with $\alpha = 1$. 

\medskip

A map $f : \mathcal{U} \to \r^m$ is then said to be  of class $C^{k,\alpha}$, where $0 < \alpha <1$ and $k\in \mathbb{N}$
if  $f \in C^{k}(\mathcal{U}, \r^m)$ and all partial derivatives of order $k$ are $\alpha$-Hölder continuous. This is a Banach space
for the norm
$$
    \|f\|_{C^{k,\alpha}(\mathcal{U})} = \|f\|_{C^{k}(\mathcal{U})}  +  \max_{\beta}  [D^{\beta}f]_{C^{0,\alpha}(\mathcal{U})}, 
$$

where $\beta\in \mathbb{N}^n$ runs through all multi-indices of order $k$.
Concerning these spaces we will need  the following results:

\begin{proposition}\label{prop.holder}
 Let \  $\mathcal{U} \subset \r^n$ be a bounded domain with Lipschitz boundary and fix 
 $k\in \mathbb{N} \cup \{0\}$ and $0\leq \alpha \leq 1$. 
 
 \smallskip
 
 (A) If $f_1, f_2 \in C^{k,\alpha}(\mathcal{U}, \r)$, then we also
 have $f_1f_2\in C^{k,\alpha}(\mathcal{U}, \r)$ and 
 $$
      \|f_1 f_2\|_{C^{k,\alpha}(\mathcal{U})} \leq 
      C \left(\|f_1 \|_{C^{0}(\mathcal{U})}\|f_2\|_{C^{k,\alpha}(\mathcal{U})} +
      \|f_2 \|_{C^{0}(\mathcal{U})}\|f_1\|_{C^{k,\alpha}(\mathcal{U})}
       \right )
 $$
 for some constant $C = C(\mathcal{U}, k)$.
 
  \smallskip
  
(B) If  $f \in C^{k,\alpha}(\mathcal{U}, \r)$ and $f \geq a > 0$ in $\mathcal{U}$, then
 $1/f\in C^{k,\alpha}(\mathcal{U}, \r)$ and
 $$
       \left\|\frac{1}{f}\right\|_{C^{k,\alpha}(\mathcal{U})} \leq   C \cdot a^{-(k+2)} 
       \|f\|_{C^{0}(\mathcal{U})}\|f\|_{C^{k,\alpha}(\mathcal{U})}
 $$
  for some constant $C = C(\mathcal{U}, k)$
\end{proposition}

The proof is given in \cite[Theorems 16.28 and 16.29]{CDK}.

\medskip
 
 \begin{proposition}\label{prop.intholder}
Let $\mathcal{U} \subset \r^n$ be a an open set and $(S,\mu)$ be an arbitrary measure space.
Let $h : \mathcal{U}\times S \to \r$ be  a measurable function such that for a.e. $s \in S$ we have
$$
   \|h_s\|_{C^{k,\alpha}(\mathcal{U})}  \leq  m(s),
$$
where $h_s(x) = h(x,s)$ and $m \in L^1(S, d\mu)$, then the function $H : \mathcal{U} \to \r$
defined by
$$
  H(x) = \int_S h(x,s) ds
$$
belongs to    ${C^{k,\alpha}(\mathcal{U})}$ and satisfies 
$ \|H\|_{C^{k,\alpha}(\mathcal{U})}  \leq    \|m\|_{L^1(S)}$.
\end{proposition}
 
\textbf{Proof.}    Let us first consider the case $k=0$, we have for any pair of
distinct points $x$ and $y$ in $\mathcal{U}$
$$
 \frac{|H(y) - H(x)|}{|y-x|^\alpha}  \leq \int_S  \frac{|h((y,s) - h(x,s)|}{|y-x|^\alpha} d\mu (s) .
$$
On the other hand, using Lebesgue's dominated convergence theorem, we easily prove that
$$
     \|H\|_{C^{k}(\mathcal{U})} \leq \int_S    \|h_s\|_{C^{k}(\mathcal{U})} d\mu (s).
$$
From the two inequalities above, we conclude that 
$$
     \|H\|_{C^{k,\alpha}(\mathcal{U})} \leq \int_S    \|h_s\|_{C^{k,\alpha}(\mathcal{U})} d\mu (s) \leq  
        \|m\|_{L^1(S)}.
$$
\qed

 \medskip 
 
We also have the following result on composition (see Theorem 16.31  in \cite{CDK}).

\begin{theorem}
Let $\mathcal{U} \subset \r^n$ and  $\mathcal{V} \subset \r^m$ be bounded open sets with Lipschitz boundaries, and let
 $h \in C^{k,\alpha}(\mathcal{V}, \r^p)$ and $f \in C^{k,\alpha}(\mathcal{U},\mathcal{V})$.
Assume also that $f$ is Lipschitz continuous, then
$h \circ f \in C^{k,\alpha}(\mathcal{U}, \r^p)$.
\end{theorem}

\medskip

\section{Regularity of the Binet-Legendre metric} \label{rec.reg}  

Theorem  A    would be an empty shell without concrete criteria
implying the $C^{k, \alpha}$ regularity of the Binet-Legendre metric. The goal of this section
is precisely to provide such criteria. We first proof Theorem B which is a result for 
$C^{k, \alpha}$ Finsler metric and then we discuss a generalized regularity condition,
called  $(k, \alpha)$-partial smoothness, for Finsler metric and we prove
a generalization of Theorem B for this class of metrics.

\medskip
 
\textbf{Proof of Theorem B.}    We can assume without loss of generality that $\mathcal{U} \subset \r^n$ is
a bounded convex domain and $F : \mathcal{U} \times \r^n \to \r$ is a continuous Finsler structure
such that 
\begin{equation}\label{Fbilip}
  \frac{1}{C_0}|v | \leq F(x,v) \leq  C_0 |v|
\end{equation}
for some constant $C_0>0$ and  all $(x, v) \in \mathcal{U}\times \r^n$.  Using Formula (\ref{eq.BL2})  and    Proposition \ref{prop.holder}, we conclude from that  
$x \mapsto  g^{i j}(x)$ is a function of class  $C^{k, \alpha}$.
The condition  (\ref{Fbilip}) on $F$ implies also that $\det(g^{ij})$ is bounded below away from zero (see e.g. \cite[Proposition 12.1]{MT2012}).
Using Proposition \ref{prop.holder} again, we then conclude that 
the inverse matrix $g_{ij}(x)$  is also a   $C^{k, \alpha}$ function of $x$ (see also Corollary 16.30 in \cite{CDK}). 

\qed

\bigskip

We will generalize Theorem B to more general Finsler metrics. To this aim we need the following

\medskip

\textbf{Definition.} Let $k\in \mathbb{N} \cup \{0\}$  and $\alpha\in  [0, 1]$. 
Let  $F$ be a continuous Finsler metric on the domain $\mathcal{U} \subset \r^n$. We say that the metric $F$ 
is $(k,\alpha)$-\emph{partially smooth}, if    there exists a   map $A : \mathcal{U} \times \r^n \to \r^n$  such that 
\begin{enumerate}[{\rm i.)}]
  \item   $A(x,\lambda v) = \lambda A(x,v)$ for all $(x,v) \in \mathcal{U} \times \r^n$ and all $\lambda\ge 0$; 
  \item  For any $x \in \mathcal{U}$, the map $A_x : \r^n \to \r^n$ defined by $v \mapsto A(x,v)$ is bilipschitz; 
  \item   There exists a constant $C > 0$, such that 
\begin{equation}\label{cond.partsmooth}
     \|h_u\|_{C^{k,\alpha}(\mathcal{U})} +  \|A_u\|_{C^{k,\alpha}(\mathcal{U})} +  \|J_u\|_{C^{k,\alpha}(\mathcal{U})} 
     \leq C
\end{equation}
  for all $u \in S^{n-1}$, where\footnote{we denote by $A_u(x) = A(x,u)$ and by $A_x(u) = A(x,u)$, this should not pose
  any problem.}
   $A_u(x) = A(x,u)$,  $J_u(x) = J(x,u) =  |\det (DA_x) (u)|$ is the Jacobian of the map $A_x$
  and $h_u(x) = F(x,A(x,u))$.  
\end{enumerate}  

\medskip

Note that this  notion of  partial smoothness is slightly different from (but strongly related to)  
the notion of $C^k$-partial-smoothness introduced in \cite[Definition 2.1.]{MT2012}.  

\smallskip

Observe also  that if we assume the maps $A_x$ to be the identity, that is $A(x,v) = v$,
then this condition reduces to the standard $C^{k, \alpha}$ regularity. 
In the general case, the  role of the map $A(x,v)$ is to allow a Finsler metric 
to be regular in $x$ after  twisting  the parameter $v$, see
\cite[\S 2]{MT2012} for some explanations and explicit examples.

\medskip

\begin{theoremBB}\label{reg.BL}
Let $\mathcal{U}$ be a bounded convex domain in $\r^n$ and let $F : \mathcal{U}\times \r^n \to \r$
be a $(k, \alpha)$ partially smooth Finsler metric in $\mathcal{U}$ satisfying the 
inequalities (\ref{Fbilip}) 
for some constant $C_0$. Then the Binet-Legendre metric  of $F$  is of class $C^{k, \alpha}$.
\end{theoremBB}

Note that this result  generalizes  Theorem B.

\bigskip

\textbf{Proof.}  The argument will be based on Proposition \ref{prop.holder},  together with a 
twisted version of  (\ref{eq.BL2}). Let us set for any $x\in \mathcal{U}$
$$
  \Omega'_x = A^{-1}_x(\Omega_x) = \{v' \in \r^n \mid F(x, A_x(v')) < 1 \}.
$$
We then have with $J(x,v') = |\det (DA)_x (v')|$
\begin{align*}
 \lambda({\Omega_x})  &= \int_{\Omega'_x} J(x,v') dv'
 \\ &= \int_{S^{n-1}}\left(\int_0^{1/F(x,A(u))}J(x,ru)r^{n-1}dr\right) du
 \\ &= \frac{1}{n} \int_{S^{n-1}}\frac{J(x,u)}{F(x,A_x(u))^{n-1}} d\sigma (u),
\end{align*}
where we have used polar coordinates $v' = ru$, with $u\in S^{n-1}$.
We have also used the fact that $A_x$ is homogenous of degree $1$
and therefore $J(x,ru)  = J(x,u)$.

\medskip

The above identity controls the denominator in the definition of
the Binet-Legendre metric. To obtain a similar formula for the
numerator, we denote by 
$$
 a_i(x,v') = (A(x,v'))_i
$$
the $i^{th.}$ coordinate of $A(x,v')$. We then have
\begin{align*}
 \int_{\Omega_x} v_iv_j dv  &= \int_{\Omega'_x} a_i(x,v')  a_j(x,v')  J(x,v') dv'
 \\ &= \int_{S^{n-1}}\left(\int_0^{1/F(x,A(u))}a_i(x,u)  a_j(x,u)  J(x,ru)r^{n+1}dr\right) du
 \\ &= \frac{1}{n+2} \int_{S^{n-1}}\frac{a_i(x,u)  a_j(x,u)  J(x,u)}{F(x,A_x(u))^{n-1}} d\sigma (u),
\end{align*}
where we have again used  $v' = ru$ and the  homogeneity $a_i(x,ru) = r a_i(x,u)$. 
Using Propositions \ref{prop.holder} and \ref{prop.intholder} together with
the inequalities (\ref{Fbilip}) and (\ref{cond.partsmooth}), we  conclude as in the  proof of Theorem B.

\qed

\medskip

We then clearly have the following

\begin{corollaryCC} 
A  distance preserving bijective map  between two  $(k,\alpha)$-partially smooth  Finsler manifolds
with $k+\alpha>0$ is a  $C_{loc}^{k+1,\alpha}$-diffeomorphism.
\end{corollaryCC}
 
 \medskip

\section{Some applications of our results}  \label{secAppl}

\subsection{Finsler manifolds with non trivial dilation} 

A direct  consequence of Corollary C  is the following
\begin{corollary}
Let $(M,F)$ be a $C^{0,\alpha}$ Finsler manifold with $\alpha >0$. Assume that $F$ is forward complete
and that there exists a map $\phi : M \to M$ such that $d_F(\phi(x),\phi(y)) = a d_F(x,y)$ for any points
$x,y \in M$ and some constant $a \neq 1$. Then $(M,F)$  is isometric to a $\r^n$ with a  Minkowski metric and 
the isometry from $M$ to $\r^n$ is a diffeomorphism.
\end{corollary} 

\textbf{Proof.}  The map $\phi$ is a distance preserving map from $(M,aF)$ to $(M, F)$. Since $F$
is H\"older continuous, one concludes from Corollary C  that $\phi$ is of class $C^{1,\alpha}$ and
$\phi^*F = a F$. The results follows then from Theorem 6.1 in \cite{MT2015}.

\qed

\subsection{On the group of isometries}

Let $(M,F)$ be a $C^0$ connected  Finsler manifold and 
 $\mathrm{Isom}(M,d_F)$ the group of its  distance preserving bijective maps $\phi : M \to M$. We then have

\begin{proposition}
 The group  $\mathrm{Isom}(M,d_F)$ is  locally compact    for the compact open topology.
\end{proposition}

\textbf{Proof.} If the map  $\phi : M \to M$ preserves the Finsler distance $d_F$, then it also preserves
its symmetrization  $\rho (x,y) = \frac{1}{2}(d_F(x,y)+d_F(y,x))$. By a theorem of  Dantzig and van der Waerden,
we know that  $\mathrm{Isom}(M,\rho)$ is locally compact for the compact open topology,  see e.g.   \cite[Theorem 4.7 in Chapter 1]{KobayshiNomizu}. The result follows since    $\mathrm{Isom}(M,d_F)$ is 
clearly a closed subgroup of $\mathrm{Isom}(M,\rho)$. 

\qed

\medskip

We now conclude that the isometry group of  a mildly regular Finsler manifold is a Lie group.

\medskip

\begin{theorem}\label{isometrics}
 Let $(M,F)$ be a smooth connected $C^0$ Finsler manifold. Assume that its Binet-Legendre metric
 is of class   $C^{0,\alpha}$  with $\alpha >0$.
 Then the group of isometries of  $(M,d_F)$  is a Lie group for   the compact open topology, and its action on $M$ is by diffeomorphisms.
\end{theorem}

\textbf{Proof.}    We just proved that   $\mathrm{Isom}(M,d_F)$ is locally
compact for the compact open topology. Since every element $\phi \in \mathrm{Isom}(M,d_F)$ is of class $C^1$, we  conclude 
 from a classical theorem of  Montgomery and his collaborators, see   \cite{BochnerMontgomery}
and  \cite[Chapter 5]{MontgomeryZippin},  that $\mathrm{Isom}(M,d_F)$ is a Lie group. By Corollary C, it acts by diffeomorphisms.

\qed

\medskip

\textbf{Remarks.}

\begin{enumerate}[a.)]
  \item Actually,   the isometry group of a Finsler manifold is a Lie group even if  the Finsler metric  is  only
   of class   $C^{0}$ (and without assuming any regularity of the Binet-Legendre metric).
  This follows from the proof of  the \emph{Hilbert-Smith conjecture} for Lipschitz homeomorphisms given in
   \cite{Repovs}. The general Hilbert-Smith conjecture states  that  \emph{a locally compact topological group $G$ that  acts faithfully on 
   a connected manifold $M$ is  a Lie group},  Recent references  are \cite{Maleshich,Michael2007,MontgomeryZippin}. 
   The conjecture has been recently proved in dimension 3 in \cite{Pardon}.
   
   \item Another particular case of the Hilbert-Smith conjecture is due to  G. Martin \cite{Martin}. It  states that  
   \emph{a locally compact group acting effectively and quasiconformally on a Riemannian manifold is a Lie group}.
   Using the Binet-Legendre metric, this results immediately extends to the case of Finsler manifolds.
 
  \item If we assume the Finsler metric $F$ to be $C^2$, then the   Binet-Legendre metric is also $C^2$ and Theorem \ref{isometrics} follows from \cite[Theorem 10]{MyersSteenrod1939}.
  
  \item In the case when the  Finsler metric $F$  is $C^2$ and strongly convex,  
  Theorem \ref{isometrics}  was proved in  \cite{DengHou2002} by S. Deng  and Z.  Hou, see also \cite[\S 3.2]{Deng2012}.
\end{enumerate}

\section{A brief history of Myers-Steenrod's Theorem} \label{secHistory}

The history behind the results discussed in the present paper is somewhat
intricate. The  original 1939 paper of Myers and Steenrod clearly sets the
problem. This paper claims that (I)   \emph{every distance preserving
homeomorphism of a $C^1$Riemannian manifold $( M, g)$ is a
$C^1$ diffeomorphism} and that (II)   \emph{any closed group of isometries
of a Riemannian manifold of class $C^2$ is a Lie group.} Both theorems are
often referred to as \emph{the Myers-Steenrod Theorem}.

\medskip

The key idea of Myers and Steenrod to prove the first result  is to represent
tangent vectors in $M$ as velocity vectors of geodesics and argue that since
geodesics are a metric notion, distance preserving maps send geodesics to
geodesics and therefore act continuously on the unit tangent bundle of $M$, from
which it follows that they are of class $C^1$. The argument is correct and is
written with some details in a more modern language in \cite[Chap.1, Theorem 11.1]{Helgason}, 
but it requires the Riemannian metric to be
at least of class $C^2$.  See also  \cite[Theorem 9.1]{Petersen} for an alternative proof,
also for $C^2$ metrics,  based on systems of coordinates defined by distance
functions. 
A  readable  proof of statement (II) is given in \cite[Chap. 6, Theorem 3.4]{KobayshiNomizu}. 

\smallskip

In 1970, Calabi and Hartman pointed at the necessity in Myers-Steenrod's argument to prove the existence of
tangent directions for the geodesics of Riemannian metrics of low regularity
and they prove a version of (I) for H{\"o}lder continuous metrics.

In 1978 however, Yu. Reshetnyak in \cite{reshetnyak}  found a mistake in a regularity result for
geodesics used by Calabi and Hartman in  \cite{Calabi1970} (in the proof of Theorem 3.1 of that paper to be precise).
This  mistake was independently discovered by A. Lytchak and A. Yaman  in \cite{Lytchak2006}, who also  constructed a 
counterexample to a  statement  Calabi and Hartman used in their proof. 
In fact the result that can be  proved by methods of  Calabi and Hartman is that  \emph{a distance preserving
homeomorphism of a $C^{0,\alpha}$ Riemannian manifold $( 0 < \alpha < 1)$ is a 
$C^{1, \frac{\alpha}{2}}$-diffeomorphism},  see \cite{sabitov1993,sabitov2008}  and \cite{Lytchak2006} for more details. 

This result was improved in Reshetnyak \cite{reshetnyak} and S. Z. Shefel' \cite {shefel}, who approached the problem via the theory of quasiconformal maps. In particular,   Theorem \ref{th.regRiem}  follows from \cite[Theorem 2]{reshetnyak}   for $k=0, \ 0 <\alpha<1$, and from  \cite{shefel} for  $k\ge 1, \ 0 <\alpha<1$. The remaining cases, when the metrics have regularity $C^{0,1}$ or $C^{1,0}$ was proved in \cite{sabitov1993}.   

\smallskip 
Another approach to this problem is due to  M. Taylor who proved   Theorem \ref{th.regRiem}  for all $k\ge 0$ and  $0\le \alpha<1$. His proof, which is elegant and self-contained, is based on harmonic coordinates and the regularity theory for elliptic PDEs
(still assuming $k+\alpha > 0$).  We refer to the book \cite{sabitov2008} for more information and further developments.

\smallskip

The generalization of Myers-Steenrod theorems to Finsler manifolds appears as a  natural question, and,   
for $C^\infty$ smooth strongly  convex metrics on manifolds of dimension at least 3,  it follows from the result of F. Brickell in   \cite{Brickell}.
The result in all dimensions has then been proved in 2002 by S. Deng  \&  Z.  Hou \cite{DengHou2002}.
See also the paper B. Aradi  \& D. Kertész \cite{Aradi} for a different point of view. All these works assume the Finsler structure $F$ to be  strongly convex, which we do not do in the present paper.

\smallskip

On the other hand, the paper  \cite{Lytchak2006} by Lytchak \& Yaman also 
considers the case of H{\"o}lder
continuous Finsler manifolds, however they need a special convexity condition on the
Finsler unit ball that they call \emph{uniform convexity of type $p$}, see
\cite[definition 2.3]{Lytchak2006}. They also prove  that \emph{a distance
preserving homeomorphism of a $C^{\alpha}$ Finsler manifold $( 0 < \alpha <
1)$ of type $p$ is a $C^{1, \beta}$ diffeomorphism with $\beta =
\frac{\alpha}{p}$.} See the remark following Theorem 1.3 in  \cite{Lytchak2006}. Our main
Theorems  extend all those previous cases and are  optimal.


\end{document}